\documentclass[10pt]{article} \usepackage{amstext,amsfonts,graphics}
\newtheorem{thm}{Theorem} \newtheorem{lemma}[thm]{Lemma}
 
\newenvironment{pf} {\noindent{\sc Proof. }}{{\hfill
$\Box$}\par\vskip2\parsep} 
\newenvironment{pfof}[1]
{\par\vskip2\parsep\noindent{\sc Proof of\ #1. }}{{\hfill $\Box$}
  \par\vskip2\parsep}

\newcommand{\py}{\mbox{\bf P}\,}
\newcommand{\ey}{\mbox{\bf E}\,}
\newcommand{\rd}{{\mathbb R}^d}
 
\newcommand{\dof}{\bf\boldmath}

\newcounter{mycount}
\newenvironment{mylist}{\begin{list}{(\roman{mycount})}%
{\usecounter{mycount}\itemsep 0pt}}{\end{list}}

\title{Trees and Matchings from Point Processes}

\author{Alexander E. Holroyd\thanks{Department of Mathematics, UC
    Berkeley, CA 94720, USA. {\tt holroyd@math.berkeley.edu}. 
Research funded in part by NSF Grant DMS--0072398.}
\and Yuval Peres\thanks{Departments of Statistics and Mathematics, UC
    Berkeley, CA 94720, USA. {\tt peres@stat.berkeley.edu}.
Research funded in part by NSF Grant DMS-0104073 and a Miller
  Professorship at UC Berkeley.}}
\begin{document}

\maketitle
\renewcommand{\thefootnote}{}
\renewcommand{\thefootnote}{\arabic{footnote}} 
\begin{abstract}
A {\dof factor graph} of a point process is a graph whose vertices are the points of the process, and which is constructed from the process in a deterministic isometry-invariant way.  We prove that the $d$-dimensional Poisson process has a one-ended tree as a factor graph.  This implies that the Poisson points can be given an ordering isomorphic to the usual ordering of the integers in a deterministic isometry-invariant way.  For $d\geq 4$ our result answers a question posed by Ferrari, Landim and Thorisson \cite{ferrari-landim-thorisson}.
We prove also that any isometry-invariant ergodic point process of finite intensity in Euclidean or hyperbolic space has a perfect matching as a factor graph provided all the inter-point distances are distinct.  
\end{abstract}

\section{Introduction}
\label{intro}

Let $M$ be an isometry-invariant point process on ${\mathbb R}^d$, viewed as a random Borel measure.  We assume throughout that all point processes are simple and of finite intensity.  The {\dof support} of $M$ is $[M]=\{x\in\rd :M(\{x\})=1\}$, and (random) elements of $[M]$ are called $M$-points. 
By a {\dof factor graph} of $M$ we mean a random (directed or undirected) graph $G$ whose vertex set equals $[M]$, such that $G$ is a deterministic function of $M$, and such that the joint distribution of $M$ and $G$ is invariant under isometries of ${\mathbb R}^d$.  (We give a more formal definition at the end of the introduction).

A graph is {\dof locally finite} if no vertex has infinite degree.  A graph is a {\dof tree} if it is connected and has no cycles.  The number of {\dof ends} of a tree is the number of distinct singly infinite self-avoiding paths from any one vertex.
A {\dof directed doubly infinite path} is a directed graph isomorphic to the graph with vertex set ${\mathbb Z}$ and a directed edge from $n$ to $n+1$ for each $n$.

\begin{thm}
\label{tree}
Let $M$ be a Poisson point process on $\rd$. 
\begin{mylist}
\item
$M$ has a factor graph which is almost surely a locally finite one-ended tree.
\item
$M$ has a factor graph which is almost surely a directed doubly infinite path.
\end{mylist}
\end{thm}
\begin{figure}
\centering \resizebox{!}{3in}{\includegraphics{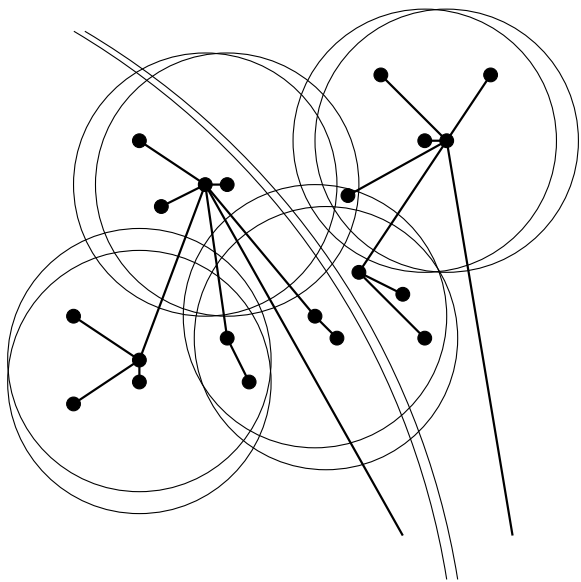}}
\caption{An illustration of the proof of Theorem \ref{tree}: seeds, cutters and blobs.}
\label{blobs}
\end{figure}

Let $\|\cdot\|$ denote the Euclidean norm, and let $B(r)=\{x\in\rd:\|x\|<r\}$ be the ball of radius $r$.
A point process $M$ is said to be {\dof non-equidistant} if there do not exist $M$-points $w,x,y,z$ with $\{w,x\}\neq\{y,z\}$ and $\|w-x\|=\|y-z\|>0$.
A graph is called a {\dof matching} if every vertex has degree $0$ or $1$, and a {\dof perfect matching} if every vertex has degree $1$.

\begin{thm}
\label{match}
Let $M$ be a non-equidistant point process in ${\mathbb R}^d$ which is invariant and ergodic under isometries.
\begin{mylist}
\item
$M$ has a factor graph in which almost surely each component is a locally finite one-ended tree.
\item
$M$ has a factor graph in which almost surely each component is a directed doubly infinite path.
\item
$M$ has a factor graph which is almost surely a perfect matching.
\end{mylist}
\end{thm}

Ferrari, Landim and Thorisson \cite{ferrari-landim-thorisson} proved by a different method that the $d$-dimensional Poisson process has a one-ended tree as a {\em translation-invariant} factor graph for $d\leq 3$, and asked whether this holds for $d\geq 4$.  Theorem \ref{tree} (i) establishes the stronger isometry-invariant statement for all $d$.  We will deduce the (ii) part of Theorem \ref{tree} from the (i) part.  This implication was also noted in \cite{ferrari-landim-thorisson}.  Theorem \ref{tree} (ii) is clearly equivalent to the assertion that the $M$-points can be given an ordering isomorphic to the usual ordering of ${\mathbb Z}$ in a deterministic isometry-invariant way. 
Such orderings have connections with the notion of point-stationarity, while the construction of perfect matchings as in Theorem \ref{match} (iii) has connections with Palm processes.  See \cite{ferrari-landim-thorisson},\cite{thorisson-book} for more details.

The main novelty of Theorem \ref{match} lies in the generality of the point process.  For the Poisson process, Theorem \ref{match} may be proved by relatively simple constructions, including one which we discuss in Section \ref{secgen}.

The minimal spanning forest is a natural factor graph.  However, in general it is unknown how many components it has, and how many ends the components have.  (Partial answers are provided in \cite{aldous-steele},\cite{alexander}).

Not every isometry-invariant ergodic point process on ${\mathbb R}^d$ has a tree as a graph factor.  For example, consider the point set obtained by applying a uniform random translation and a uniform random rotation to ${\mathbb Z}^d$.  This process has no perfect matching as a factor graph, and no factor graph in which every component is an infinite tree if $d\geq 2$.

With additional work, the condition of non-equidistance in Theorem \ref{match} may be relaxed to the condition that the symmetry group of $[M]$ is almost surely trivial.  The idea is to define an isometry-invariant ``index function'' on $M$-points, and use it to break ties between distances.  Indeed, these ideas may used to obtain a necessary and sufficient condition for an ergodic process $M$ on $\rd$ to have a perfect matching as a factor graph (for example).  Let $\Sigma$ be the (random) symmetry group of $[M]$, and consider the quotient point process (on a random manifold) $M/\Sigma$.  Then it may be shown that $M$ has a perfect matching as a factor graph if and only if the support of $M/\Sigma$ is almost surely of even or infinite cardinality.

Theorems \ref{tree} extends to other amenable spaces in place of $\rd$, with the same proof.  In non-amenable spaces, no invariant one-ended tree exists \cite{adams-lyons},\cite{blps}, and in certain spaces no invariant tree exists \cite{adams-spatzier},\cite{pemantle-peres}.  Theorem \ref{match} extends to other spaces as follows.
\begin{thm}
\label{general}
Let $\Lambda$ be a locally compact metric space and let $\Gamma$ be a transitive unimodular group of isometries of $\Lambda$.  Let $\nu$ be a $\Gamma$-invariant Borel measure on $\Lambda$ which is finite on bounded sets.  Let $M$ be a non-equidistant point process of finite intensity with respect to $\nu$, which is invariant and ergodic under $\Gamma$.  Then statements (i)--(iii) of Theorem \ref{match} hold.
\end{thm}
The proof of Theorem \ref{match} extends directly to Theorem \ref{general} in the case of amenable spaces, but breaks down for non-amenable spaces such as hyperbolic space.  We therefore give a different argument which is valid in amenable and non-amenable settings.  The argument is based on the simple and appealing idea of iteratively matching mutually nearest neighbors, and is of interest even in the case of $\rd$.

Theorems \ref{tree},\ref{match},\ref{general} are proved Sections \ref{sectree},\ref{secmatch},\ref{secgen} respectively.

\vspace{\parsep}
{\noindent\bf Formal definition of a factor graph.}  A {\dof factor graph} of a point process $M$ on $\rd$ is a mapping, commuting with isometries, that assigns to a countable set $V\subseteq\rd$ a graph $G=G(V)$ with vertex set $V$, defined for almost all $V$ with respect to the law of $[M]$.  We do not need to impose a $\sigma$-algebra on the space of graphs; the measurability with respect to the law of $M$ of the events in Theorem \ref{tree}--\ref{general} will be evident from the proofs.  Similar remarks apply to clumpings as defined in Section \ref{sectree}.

\section{Poisson Trees}
\label{sectree}

\begin{lemma}
\label{locally-finite}
Let $M$ be an isometry-invariant point process on ${\mathbb R}^d$, and let $G$ be a factor graph all of whose components are trees with at most $r$ ends for some constant $r<\infty$.  Then $G$ is locally finite almost surely.
\end{lemma}

In this section we will need only the special case $r=1$. 
The proof of Lemma \ref{locally-finite} uses a version of the ``mass transport principle'', Lemma \ref{mass-transport}.  See \cite{blps},\cite{hyperbolic},\cite{haggstrom-mt} for background.
A {\dof mass transport} is a non-negative measurable function $T(x,y,M)$ on $\rd\times\rd\times{\cal M}$ (where ${\cal M}$ is the space of Borel measures on $\rd$) which is non-zero only when $x,y$ are $M$-points, and which is isometry-invariant in the sense that $T(\theta x,\theta y,\theta M)=T(x,y,M)$ for any isometry $\theta$ of $\rd$.  We think of $T(x,y,M)$ as the mass sent from $x$ to $y$ when the point configuration is $M$.  For Borel sets $A,B\subseteq\rd$ we write $$t(A,B)=\ey \sum_{\begin{array}{c}\\[-.2in]\scriptstyle x\in A\cap[M],\\[-.04in] \scriptstyle y\in B\cap[M]\end{array}}
 T(x,y,M)$$ for the expected total mass sent from $A$ to $B$, and we write $K=[0,1)^d\subseteq\rd$ for the unit cube.

\begin{lemma}
\label{mass-transport}
Let $M$ be an isometry-invariant point process on ${\mathbb R}^d$, and let $t$ be a mass transport.  We have
$$t(K,\rd)=t(\rd,K).$$
\end{lemma}

\begin{pf}
The isometry-invariance of $T$ implies isometry-invariance of $t$.  Hence, using monotone convergence to exchange expectations and sums, we have
$$
t(K,\rd)=\sum_{z\in{\mathbb Z}^d}t(K,K+z)=\sum_{z\in{\mathbb Z}^d}t(K-z,K)
=t(\rd,K).$$
\vspace*{-1\in}
\end{pf}

\begin{pfof}{Lemma \ref{locally-finite}}
Consider the mass transport in which $T(x,y,M)=1$ whenever $G$ has a singly infinite self-avoiding path from $x$ which includes the edge $(x,y)$, and $T(x,y,M)=0$ otherwise.  Thus each vertex sends one unit of mass to each of at most $r$ of its neighbors, while for every edge $(x,y)$ in an infinite component, either $x$ sends one unit to $y$ or $y$ sends one unit to $x$, or both. 
Applying Lemma \ref{mass-transport}, the assumption of finite intensity implies that the expected total mass received by all $M$-points in $K$ is finite, so in particular it follows that no vertex can have infinite degree.
\end{pfof}

A {\dof clumping} of a point process $M$ is a sequence ${\cal P}_1,{\cal P}_2,\ldots$ of successively coarser partitions of $[M]$, defined from $M$ in a deterministic isometry-invariant way.  (More precisely, a clumping is a mapping, commuting with isometries, that assigns to a countable set $V\subseteq\rd$ a sequence of successively coarser partitions of $V$, defined for almost all $V$ with respect to the law of $[M]$).  We call the elements of the partitions {\dof clumps}, and we call the clumping {\dof locally finite} if all clumps are finite.
A {\dof component} of a clumping is subset of $[M]$ which is the limit of some increasing sequence of clumps $A_1,A_2,\ldots$, where $A_k\in{\cal P}_k$.
A clumping is {\dof connected} if it has only one component.

\begin{lemma}
\label{clumping}
Let $M$ be a non-equidistant isometry-invariant point process on $\rd$.  If $M$ has a clumping which is almost surely connected and locally finite then it has a factor graph which is a almost surely a locally finite one-ended tree.
\end{lemma}

\begin{pf}
Let ${\cal P}_1,{\cal P}_2,\ldots$ be a connected locally finite clumping.  We will define for each clump $A$ a distinguished element $x\in A$ called the {\dof leader} of $A$, and we will construct a factor graph $G$.  We do this inductively as follows.
First consider a clump $A\in {\cal P}_1$.  Choose the leader $x$ of $A$ as follows.  If $|A|=1$, let $x$ be the unique element.  Otherwise, let $x',x''$ be the unique pair of $M$-points in $A$ whose Euclidean distance is minimum.  Then let $x$ be the one of $x',x''$ which minimizes $\min\{\|x-y\|:y\in[M]\setminus\{x',x''\}\}$.  The fact that $M$ is a non-equidistant point process ensures that all the minima involved are unique almost surely.  Let $G$ have an edge from $x$ to each of the other elements of $A$.
Apply the same construction to every $A\in{\cal P}_1$.

Now suppose that leaders have been defined for all clumps in the partitions ${\cal P}_1,\ldots,{\cal P}_{k-1}$.  Consider a clump $A\in {\cal P}_k$, let $B_1,\ldots,B_m$ be the clumps of ${\cal P}_{k-1}$ which are subsets of $A$, and let $y_1,\ldots,y_k$ be their respective leaders.  Choose the leader $x$ of $A$ from among $y_1,\ldots,y_k$ in a deterministic isometry-invariant way by the same procedure as above.  Let $G$ have an edge from $x$ to each $y_i\neq x$.
Apply the same construction to every $A\in{\cal P}_k$.

 By the construction $G$ is clearly a factor graph.
  Also, if $A$ is any clump, then the subgraph of $G$ induced by $A$ is a finite tree, and only the leader of $A$ has any edges to vertices outside $A$.  This implies immediately that all components of $G$ are one-ended trees or finite trees.  But $G$ is connected since the clumping is connected, hence $G$ is a one-ended tree.  Finally $G$ is locally finite by Lemma \ref{locally-finite}.
\end{pf}

\begin{lemma}
\label{one-two}
Let $M$ be a non-equidistant isometry-invariant point process on $\rd$.  If $M$ has a factor graph which is almost surely a locally finite one-ended tree then it has a factor graph which is almost surely a directed doubly infinite path.
\end{lemma}

\begin{pf}
Let $G$ be such a one-ended tree.  First
order the children of each vertex in order of distance to the parent.
Then order all the vertices $G$ according to
depth-first search; that is, each vertex precedes all its children, while if $x$ precedes its sibling $y$ then all descendants of $x$ precede all descendants of $y$.  (See \cite{ferrari-landim-thorisson} for another description of this construction).
\end{pf}

\begin{pfof}{Theorem \ref{tree}}
Suppose without loss of generality that $M$ has intensity $1$.
Throughout we will use $C_1,C_2,\ldots$ to denote constants in $(0,\infty)$ depending only on $d$.  We will construct a connected clumping of $M$; then the theorem follows immediately from Lemmas \ref{locally-finite},\ref{clumping},\ref{one-two}.

Our argument is based on a construction in \cite{blps}.  For each integer $k\geq 1$, let $a_k=\exp\left[-k\left(1-\frac{1}{2d}\right)\right]$.  Call an $M$-point $x$ a {\dof $k$-seed} if there is another $M$-point  within Euclidean distance $a_k$ of $x$.  Clearly, the $k$-seeds form an isometry-invariant point process of intensity
\begin{equation}
\lambda_k:=1\times \py[M(B(a_k))\geq 1]=1-e^{-C_1 a_k^d}\leq C_1 a_k^d = C_1 e^{-k(d-1/2)}.
\label{intensity}
\end{equation}
Now let $r_k=e^k$, and define a {\dof $k$-cutter} to be any subset of ${\mathbb R}^d$ of the form $\{y: \|y-x\|=r_k\}$, where $x$ is a $k$-seed.  Let $W_k$ be the union of all $k$-cutters, and
define a {\dof $k$-blob} to be any connected component of ${\mathbb R}^d\setminus \bigcup_{j\geq k} W_j$.  Clearly, every $M$-point lies in exactly one $k$-blob for each $k$ almost surely, and every $k$-blob is a subset of exactly one $(k+1)$-blob.  (Note that $k$-cutters typically occur in nearly-coincident pairs.  This fact makes the pictures rather odd, but neither helps nor hinders our proof.  The anomaly could be avoided at the expense of a less convenient definition of $k$-seeds).

For each $k$, define a partition ${\cal P}_k$ of $[M]$ by declaring two $M$-points to be in the same clump of ${\cal P}_k$ if they lie in the same $k$-blob.  Clearly ${\cal P}_1,{\cal P}_2,\ldots$ form a clumping; we must check that it is locally finite and connected.

First we claim that almost surely all blobs are bounded.  This will imply immediately that the clumping is locally finite.
It is sufficient to check that almost surely all blobs which intersect $B(1)$ are bounded.
Let 
$$V_k=\{B(1) \text{ is enclosed by some $k$-cutter}\}.$$
We will show that
$\py[V_k]\rightarrow 1$
as $k\rightarrow \infty$.  This implies that $V_k$ occurs for infinitely many $k$ almost surely, and this implies the claim, since any $j$-blob which intersects $B(1)$ must then be enclosed by some $k$-cutter for some $k\geq j$.  By the definition of a $k$-cutter, $V_k$ equals the event that $B(r_k-1)$ contains some $k$-seed.  Now, we have
\begin{eqnarray*}
&\py[B(a_k/2) \text{ contains some $k$-seed}]\geq \py[M(B(a_k/2))=2]& \\
&\geq  e^{-C_2 a_k^d} (C_2 a_k^d)^2/2
\geq  C_3 a_k^{2d}&
\end{eqnarray*}
Moreover, if $A,B\subseteq {\mathbb R}^d$ are sets at distance at least $2a_k$ from each other, then the events $\{A \text{ contains some $k$-seed}\},\{B \text{ contains some $k$-seed}\}$ are independent.  For $k\geq 2$ we may clearly find $\lceil C_4 (r_k/a_k)^d \rceil$ balls of radius $a_k/2$ lying in $B(r_k-1)$ and spaced at distance at least $2a_k$ from each other.
Hence
\begin{eqnarray*}
&\py[V_k]=\py[B(r_k-1) \text{ contains some $k$-seed}]\geq 1-(1-C_3 a_k^{2d})^{C_4 (r_k/a_k)^d} & \\
& \geq  1- e^{-C_5 a_k^d r_k^d}
 = 1- e^{-C_5 e^{k/2}} \rightarrow 1,&
\end{eqnarray*}
establishing the above claim.

To prove that the clumping is connected, it suffices to prove that for every fixed $\ell>0$, almost surely all $M$-points in $B(\ell)$ lie in the same clump of ${\cal P}_k$ for some $k$.  By the construction it is enough to show that all such $M$-points lie in a single blob, and this in turn follows if almost surely all of $B(\ell)$ lies in some blob.
Let 
$$U_k=\{B(\ell) \text{ intersects some $k$-cutter}\}.$$
We have
\begin{eqnarray*}
\sum_{k=1}^{\infty} \py[U_k] &=& \sum_{k=1}^{\infty} \py[B(r_k+\ell)\setminus B(r_k-\ell) \text{ contains some $k$-seed}] \\
&\leq & \sum_{k=1}^{\infty} \lambda_k C_6 \ell r_k^{d-1} 
\leq  C_7 \ell \sum_{k=1}^{\infty} e^{-k/2} \qquad\text{(by (\ref{intensity}))}\\
&<& \infty,
\end{eqnarray*}
so the Borel-Cantelli lemma implies that $\py[U_k \text{ occurs for infinitely many }k]$ $=0$, and hence $B(\ell)$ lies in some blob as required.
\end{pfof}

\section{Forests and Matchings}
\label{secmatch}

\begin{lemma}
\label{equiv} 
Let $d\geq 1$ and let $M$ be a non-equidistant point process in ${\mathbb R}^d$ which is invariant and ergodic under isometries.  The following are equivalent.
\begin{mylist}
\item
$M$ has a clumping which almost surely is locally finite and has each component infinite.
\item
$M$ has a factor graph in which almost surely each component is a locally finite one-ended tree.
\item
$M$ has a factor graph in which almost surely each component is a directed doubly infinite path.
\item
$M$ has a factor graph in which almost surely each component is a locally finite one-ended or two-ended tree.
\end{mylist}
\end{lemma}

\begin{pf}
The implications (i)$\Rightarrow$(ii) and (ii)$\Rightarrow$(iii) follow immediately by applying the proofs of Lemmas \ref{clumping},\ref{one-two} to each component.  The implication (iii)$\Rightarrow$(iv) is a triviality since a doubly infinite path is a two-ended tree.  Therefore it is sufficient to prove (iv)$\Rightarrow$(i).

Suppose $G$ is a factor graph as in (iv).  We will treat each component separately.  Let $C$ be a component of $G$.  If $C$ is a one-ended tree, let $\pi$ be any singly infinite path in $C$, and define
$$L(C)=\limsup \{\| x-y \| :(x,y) \text{ is an edge of }\pi \}.$$
Since any two singly infinite paths in $C$ must eventually coalesce, $L(C)$ does not depend on the choice of $\pi$.  On the other hand if $C$ is a two-ended tree, then it has a unique doubly infinite self-avoiding path, which we call the {\dof trunk}.  If the trunk is deleted, only finite components remain.  Define
$$L(C)=\sup\{\|x-y\|:(x,y) \text{ is an edge of the trunk}\}.$$

Let $k$ be a positive integer.  Define for each component $C$
$$L_k(C)=\left\{
\begin{array}{ll}
L(C)-k^{-1} & \text{ if }L(C)<\infty, \\
k & \text{ if }L(C)=\infty.
\end{array}
\right.
$$
Let $G_k$ be the graph obtained from $G$ by deleting every edge $(x,y)$ for which $\|x-y\|\in [L_k(C),L(C))$, where $C$ is the component containing $(x,y)$.  
Let ${\cal P}_k$ be the partition of $[M]$ induced by the components of $G_k$.  Clearly ${\cal P}_1,{\cal P}_2,\ldots$ form a clumping;
we claim that it is locally finite and has the same components as $G$; this will establish (i).

First note that there is no edge $(x,y)$ for which $\|x-y\|=L(C)$ where $C$ is the component of $G$ containing $(x,y)$.  To see this, note that by non-equidistance, there can be at most one such edge in each component of $G$.  But now consider the mass transport in which every vertex in component $C$ sends one unit of mass to each of $x,y$ if there is such an edge $(x,y)$ in $C$.  Applying Lemma \ref{mass-transport} gives a contradiction since each vertex sends out at most two units, but such vertices $x,y$ would receive infinite mass.

It follows that if $x,y$ are vertices in the same component $C$ then there exists some $k$ such that no edge in the path from $x$ to $y$ has length lying in $[L_k(C),L(C))$, so $x,y$ lie in the same clump of ${\cal P}_k$.  Hence the components of the clumping are the components of $G$.

It remains to show that the all components of $G_k$ are finite, since this will imply that the clumping is locally finite.  If $C$ is a one-ended component of $G$, the definition of $L(C)$ implies that every infinite path in $G$ has some edge with length in $[L_k(C),L(C))$, so all components of $G_k$ which lie in one-ended components of $G$ are finite.  Let $H$ be the graph consisting of all trunks of two-ended components of $G$, and let $H_k=G_k \cap H$.  Clearly all components of $H_k$ are doubly infinite paths, singly infinite paths, or finite paths.  We claim that in fact the first two possibilities can be ruled out, and this implies that all components of $G_k$ must be finite.  We prove the claim as follows.  Firstly, the definition of $L(C)$ for a two-ended component $C$ of $G$ implies that the trunk of $C$ has some edge with length in $[L_k(C),L(C))$, so components of $H_k$ cannot be doubly infinite paths.  Secondly, consider the mass transport in which every vertex in a singly infinite path component of $H_k$ sends one unit of mass to the end point of the singly infinite path.  Each vertex sends out at most one unit, but if singly infinite paths existed then their end points would 
receive infinite mass, violating Lemma \ref{mass-transport}.
\end{pf}

\begin{pfof}{Theorem \ref{match}}
The {\dof minimal spanning forest} $S$ of $M$ is the factor graph obtained from the complete graph on $[M]$ by deleting every edge which is the longest in some cycle.  It is proved in \cite{alexander} that all components of $S$ are one-ended or two-ended trees almost surely.  And $S$ is locally finite by Lemma \ref{locally-finite}.  Hence Theorem \ref{match} (i),(ii) follow from Lemma \ref{equiv}.

Let ${\cal P}_1,{\cal P}_2,\ldots$ be a clumping as in Lemma \ref{equiv} (i).  We will construct a perfect matching, establishing Theorem \ref{match} (iii).  First, consider a clump $A\in{\cal P}_1$ which has at least two $M$-points, match the two closest $M$-points of $A$, then remove them and repeat until $A$ has at most one unmatched $M$-point.  Do this for every clump $A\in{\cal P}_1$.  Next apply the same construction to each clump $A\in
{\cal P}_2$ which has at least two unmatched $M$-points, and so on indefinitely.  It is clear that each component of the clumping contains at most $M$-point that is never matched.  We claim that in fact there are no such $M$-points.  To prove this consider the mass transport in which $x$ sends one unit to $y$ if $x,y$ are in the same component of the clumping and $y$ is never matched.  Since all components of the clumping are infinite, any $M$-point never matched would receive infinite mass, contradicting Lemma \ref{mass-transport}.
\end{pfof}

\section{Iterated Nearest Neighbor Matching}
\label{secgen}

In this section we prove Theorem \ref{general}.  The mass transport principle in Lemma \ref{mass-transport} extends to the more general setting, with any bounded Borel set $K$ in place of the unit cube.  In the case of the hyperbolic plane the result is a special case of Theorem 5.2 in \cite{hyperbolic}, and the proof in \cite{hyperbolic} extends to our more general setting.
Lemma \ref{equiv} extends to the general setting with essentially the same proof.
 In the case when $\Gamma$ is amenable, the proof of Theorem \ref{match}
also extends, since all the components of any invariant random forest
must have either one or two ends.
In the non-amenable case, this argument breaks down;
although it is believed that the components of the
 (wired) minimal spanning forest all have at most two ends, this has not
been proved. Thus we use a different construction, which is of interest
even in Euclidean space.

Let $\rho$ be the metric on $\Lambda$.  By a {\dof descending chain} we mean a sequence of distinct $M$-points $x_1,x_2,\ldots$ for which the distances $\rho(x_i,x_{i+1})$ form a (strictly) decreasing sequence.  We will see that the argument becomes simpler in the case when $M$ has no descending chains.  In particular it was proved in \cite{haggstrom-meester} that Poisson processes have no descending chains almost surely.  On the other hand, there do exist isometry-invariant point processes with descending chains.  For example, start with a Poisson process on $\rd$, construct a one-ended tree as a factor graph (Theorem \ref{tree}), and then add extra points along the edges in such a way that every singly infinite path becomes a descending chain.

\begin{pfof}{Theorem \ref{general}}
We start by proving that $M$ has a perfect matching as a factor graph.  Then we will deduce the statements (i),(ii) in Theorem \ref{match}.

Call a pair of $M$-points $x,y$ {\dof mutually closest} if $x$ is the closest $M$-point to $y$ and $y$ is the closest $M$-point to $x$.  We may construct a (not necessarily perfect) matching $G$ as a factor graph of $M$ by the following procedure.
 First match all mutually closest pairs of $M$-points, and remove such points,
then match and remove all mutually closest pairs among the remaining points, and repeat indefinitely.  Let $N$ be the process of all $M$-points which are never matched by this procedure (that is, that have degree $0$ in $G$).  Clearly $N$ is an isometry-invariant ergodic point process, so in particular it has almost surely infinitely many points or almost surely no points.  The matching $G$ is perfect if and only if $N$ has no points.  (The above procedure was suggested by Dana Randall, and is noted in \cite{ferrari-landim-thorisson}.  It has the following informal interpretation.  Imagine a growing ball centered at each $M$-point, such that at time $t$ each ball has radius $t$.   Every time two balls meet, they are annihilated and their centers are matched).

Consider the directed factor graph of $H$ of $N$ in which there is a directed edge from each $N$-point to its closest $N$-point.  It is easy to see that $H$ has no cycles except of size $2$ (these being exactly the mutually closest pairs of $N$-points), and that every finite component of $H$ contains exactly one cycle of size $2$.  We claim that in fact $N$ has no mutually closest pairs, so $H$ has no finite components.  To see this, suppose that $x,y$ are mutually closest $N$-points.  This is equivalent to the statement that the set $J=\{z\in\Lambda: \rho(z,x)\leq \rho(x,y) \text{ or }\rho(z,y)\leq \rho(x,y)\}$ contains no $N$-points other than $x,y$.  Since $J$ is bounded, it contains only finitely many $M$-points almost surely, and hence at some finite (random) stage of the matching procedure above, $J$ contains no unmatched $M$-points other than $x,y$.  But then $x,y$ will be matched at the next stage, which contradicts the assumption that they are $N$-points.

The above argument shows that if $G$ is not a perfect matching then $H$ is non-empty and every component of $H$ is an infinite tree.  (Note in particular that the latter is possible only if $M$ has descending chains).  We claim that the components of $H$ are one-ended or two-ended trees.  Once this is established we obtain a perfect matching of $M$ by using Lemma \ref{equiv} (extended to $\Lambda$) to perfectly match $N$, and then combining this with $G$.

To prove the above claim, define the {\dof backbone} $B$ of $H$ to be the directed subgraph of $N$ whose edge set is the union of all (not necessarily directed) doubly infinite self-avoiding paths of $H$.  Note that all vertices incident to $B$ have degree at least $2$ in $B$.  The claim is equivalent to the assertion that $B$ has no vertices of degree greater than $2$.  Consider the mass transport in which $x$ sends one unit of mass to $y$ if $B$ has a directed edge from $x$ to $y$.  Note that the degree of a vertex in $B$ equals the total mass sent out plus the total mass received.
By the construction of $H$, each vertex sends out at most one unit.  Let $K$ be a fixed bounded Borel set, let $V$ be the set of $N$-points in $K$ incident to $B$, and let $D=\sum_{x\in V} \mbox{deg}_B(x)$.  Then the mass transport principle yields
$$\ey D\leq 2{\ey|V|}.$$
On the other hand, if $B$ has vertices of degree greater than $2$ then
$$\ey[D-2|V|]>0,$$
a contradiction, proving the claim.

The above argument shows that every $M$ satisfying the conditions of Theorem \ref{general} has a perfect matching as a factor graph.  We deduce that statements (i),(ii) of Theorem \ref{match} hold as follows.  Construct a perfect matching, then delete one $M$-point from each matched pair (choosing which one as in the choice of leaders in the proof of Lemma \ref{clumping}, for example), then construct a perfect matching as a factor graph of the process of remaining points, and repeat indefinitely.  This gives a locally finite clumping with all components infinite, so we can apply Lemma \ref{equiv} (extended to $\Lambda$).
\end{pfof}

\section{Open Problems}
\label{probs}

\begin{itemize}

\item
For what other point processes do the conclusions of Theorem \ref{tree} hold?  In particular, do they hold for every non-equidistant ergodic point process on $\rd$?
\item
Consider a perfect matching which is a factor graph of a point process.  What can be said about the lengths of the edges?  More specifically, how does the probability that some point in $B(1)$ is matched to a point outside $B(r)$ behave as $r\to\infty$?  One may ask such questions for specific matchings such as those discussed in Section \ref{secgen}, or one may ask for the optimal tail behavior over all possible perfect matchings.  The latter question is related to problems studied in \cite{steele}, for example.
\item
What other graphs are possible as factor graphs of point processes?  For example, for which $d,n$ does the Poisson process on $\rd$ have ${\mathbb Z}^n$ as a factor graph?
\end{itemize}

\section*{Acknowledgments}
We thank David Aldous, Dana Randall and Hermann Thorisson for valuable discussions.

\bibliography{phd}

\end{document}